\input amstex
\documentstyle{amsppt}
\define\tr{\text{trace}}
\define\fp{(\text{f.p.})_{z=0}}
\define\Reso{(\lambda I-A)^{-1}}
\define\reso{(\lambda-a)^{-1}}
\define\res{\text{res}}
\NoRunningHeads
\topmatter
\title Determinants of Zeroth Order Operators\endtitle
\author Leonid Friedlander and Victor Guillemin\endauthor
\affil University of Arizona\\
       Massachusetts Institute of Technology \endaffil
\endtopmatter
\document
\head 1. Introduction\endhead
In this paper we will compare two
techniques for defining regularized determinants of zeroth order
pseudodifferential operators and show that, modulo local terms,
they give the same answer. To illustrate these two techniques let
$f$ be a $C^\infty$ function on the circle with $f-1\approx 0$.
Szeg\H o proves that if $P_n$ is orthogonal projection on the
space spanned by $e^{ik\theta}$, $-n\leq k\leq n$, then for $n$
large $$\log\det P_nM_fP_n=2n\widehat{\log f}(0)+\sum k\widehat{\log
f}(k)\widehat{\log f}(-k)+O\bigl(n^{-\infty}\bigr)\tag 1.1$$ where
$M_f$ is the operator of multiplying by $f$ and $\widehat{\log f}(k)$
is the $k$-th Fourier coefficient of $\log f$. Hence by
subtracting off the "counterterm" $2n\widehat{\log f}(0)$ one gets
for the Szeg\H o-regularized determinant of $M_f$:
$$\log\det M_f =\sum k\widehat{\log
f}(k)\widehat{\log f}(-k).\tag 1.2$$

An alternative way of regularizing this determinant is by zeta
function techniques. Namely, let $Q^z:L^2(S^1)\to L^2(S^1)$ be
the operator
$$Q^ze^{in\theta}=\cases |n|^ze^{in\theta},& \text{if $n\not=0$};\\
                          0,& \text{if $n=0$}.\endcases$$
Then
$$\tr (\log M_f)Q^z=\tr M_{\log f}Q^z=2\widehat{\log
f}(0)\sum_{n=1}^\infty n^z= 2\widehat{\log f}(0)\zeta(-z).$$
Since $\zeta(z)$ is regular at $z=0$ and $\zeta(0)=-1/2$,
zeta-function regularization gives one, for the regularized
"$\log\det$" of $M_f$
$$\tr \log M_f=-\widehat{\log f}(0),\tag 1.3$$
i.e. the zeta-regularized determinant of $M_f$ is proportional to
the counterterm one had to subtract off in order to obtain the
Szeg\H o-regularized determinant of $M_f$.

This does not bode well for comparing these two methods of
regularization in more general setting; however, the right hand
sides in (1.2) an (1.3) are local expressions of the symbol of
$M_f$, and for both these methods of regularization the non-local
contributions are zero. In the paper we will show that if one
replaces $M_f$ by a zeroth order pseudodifferential operator, $B$,
then (1.2) and (1.3) are non-symbolic (i.e. non-local) functions
of $B$; however, their difference is symbolic. In other words,
modulo local terms, they give the same answer.

This is a special case of a more general result about
"Zoll operators". Let $X^d$ be a compact manifold and
$Q:C^\infty(X)\to C^\infty(X)$ a self-adjoint first order elliptic
pseudodifferential operator. $Q$ is a {\it Zoll} operator if the
bicharacteristic flow on $T^*X\setminus X$ generated by its symbol
is periodic of period $2\pi$. (To simplify the statements of some
of the results below we'll strengthen this assumption and assume
the bicharacteristic flow strictly periodic of period $2\pi$: if
the initial point of a bicharacteristic is $(x,\xi)$, the
bicharacteristic returns for the first time to $(x,\xi)$ at
$t=2\pi$.) If $Q$ is a Zoll operator, the operator
$$W=-{1\over 2\pi i}\log\exp 2\pi i Q,$$
with $0<\hbox{\rm Im}\log z\leq 2\pi$, is a zeroth order
pseudodifferential operator, and the spectrum of the operator
$Q+W$ consists of positive integers. We'll henceforth subsume
this property into the definition of "Zoll", and assume
$\hbox{\rm spec} Q={\Bbb Z}_+$. (The standard example of a
Zoll operator is the operator
$$\biggl(\Delta_{S^d}+\biggl({d-1\over
2}\biggr)^2\biggr)^{1/2}-{d-1\over 2};$$
however, there are a lot of non-standard examples as well.
See, for instance [CV].)

Let $\pi_k$ be the orthogonal projection of $L^2(X)$ onto the
$k$-th eigenspace of $Q$ and let $P_n=\pi_1+\cdots+\pi_n$.
If $B:L^2(X)\to L^2(X)$ is a zeroth order pseudodifferential
operator and $I-B$ is small then by a theorem of Guillemin and
Okikiolu [GO]
$$\log\det P_nBP_n\sim b+\sum_{k=d,
k\not=0}^{-\infty}b_kn^k+b_0\log n\tag 1.4$$
and, as above, one can define the Szeg\H o regularized determinant
of $B$ to be $e^b$. On the other hand, the expression
$$\tr \log B Q^z\tag 1.5$$
is a meromorphic function in $z$ with simple poles at $z=-d+k$,
$k=0,1,\ldots$, and one can define the zeta function
regularization of $\log\det B$ to be the finite part of this
function at $z=0$.

In section 2 we will compare these two definitions and show that,
as above, they differ by an expression that is local in $B$ and
only involves integrals of terms in the symbolic expansion of $B$
of degree $\geq -d$. Then in section 3 we will examine zeta
regularization in more detail, allowing the "regularizer" $Q$ to be
any positive definite self-adjoint first order elliptic
pseudodifferential operator (i.e.,~not necessarily a Zoll operator
as above) and prove a number of results about the "$\log\det$":
$$w_Q(B)=(\text{f.p.})_{z=0}\tr (\log B)Q^z\tag 1.6$$
for zeroth order pseudodifferential operators, $B$.
For instance we will show that the variation, $\delta w_Q$,
of this functional is local and that if $Q$ and $Q'$ are two
regularizers, $w_Q(B)-w_{Q'}(B)$ is local. (In other words,
modulo local terms, the
regularization of $\log\det B$ defined by (1.6) is independent of
the choice of $Q$.) We will also compute the multiplicative
anomaly of the regularized $\log\det B_1B_2$ defined by (1.6) and
show that it, too, is given by expressions which are local in the
symbols of $B_1$ and $B_2$.
\head 2. Szeg\H o regularized determinants\endhead
We will give a brief sketch of how (1.4) was derived in [GO]
and show how the zeroth order term in this expression is related
to (1.6). Letting $B=I-A$ the left hand side of (1.4) becomes
$$\sum_{k=1}^\infty {1\over r}\tr (P_nAP_n)^r,\tag 2.1$$
so to study the asymptotic behavior of (1.4) it suffices
to study the asymptotic behavior as $n$ tends to infinity of each of the 
summands in (2.1). To do this we will decompose the operator $A$ into
its ``Fourier coefficients'' as in the example discussed in section 1. More 
explicitely let $U(t)=\exp (itQ)$ and let
$$A_k={1\over 2\pi}\int_0^{2\pi}e^{ikt}U(-t)AU(t)dt.\tag 2.2$$
By Egorov's theorem the $A_k$'s are zeroth order pseudodifferential
operators, and the sum
$$A=\sum_{k=-\infty}^\infty A_k$$
is the ``Fourier series'' of $A$. It is shown in [GO] that this series 
converges and that the operator norms of the $A_k$'s are rapidly decreasing in 
$k$ as $k$ tends to infinity. Hence for deriving asymptotic expansions
for the summands in (2.1) we can assume that
$$A=\sum_{k=-N}^NA_k,\quad N\ \text{large}.\tag 2.3$$
Also, since $U(t)=\sum e^{int}\pi_n$,
$$A_k={1\over 2\pi}\sum_{m,n}\int_0^{2\pi} e^{ikt}e^{i(n-m)t}\pi_mA\pi_ndt
=\sum_n \pi_{n+k}A\pi_n.\tag 2.4$$
Plugging (2.3) into the $r$th summand of (2.1) and replacing each term in the 
product by the sum (2.4) one gets:
$$\tr (P_nAP_n)^r=\sum_{j_1+\cdots+j_r=0}\tr\sum_{k+\sigma(j)\leq n}
\pi_kA_{j_r}\cdots A_{j_1}\pi_k\tag 2.5$$
where $j=(j_1,\ldots,j_r)$,
$$\sigma(j)=\max(0,j_1,j_1+j_2,\ldots,j_1+\cdots+j_r),\tag 2.6$$
and the number of summands in $j$ is finite. We will use the notation
$A_j=A_{j_r}\cdots A_{j_1}$.
The asymptotics of each of the summands in (2.5) can be read off from
a theorem of Colin de Verdiere [CV] which says that
$$\tr \pi_nA_j\pi_n\sim\sum c_l(A_j)n^l.\tag 2.7$$
Moreover, Colin's theorem asserts that the terms on the right are local 
functionals of $A$ and are given explicitly by the non-abelian residues
$$c_l(A_j)=\res Q^{-(l+1)}A_j.\tag 2.8$$
Finally by plugging (2.8) into (2.7) we obtain an asymptotic expansion
$$\tr (P_nAP_n)^r\sim a_r+\sum_{k=d, k\not=0}^{-\infty}a_{r,k}n^k
+a_{r,0}\log n\tag 2.9$$
in which all terms except the constant term, $a_r$, are local functions of 
$A$.

The same argument can also be used to compute $\tr (P_nAP_n)^rQ^z$. 
Namely, by (2.5),
$$\tr (P_nAP_n)^rQ^z=\sum_{j_1+\cdots+j_r=0}\tr\sum_{k+\sigma(j)\leq n}
\pi_kA_{j_r}\cdots A_{j_1}\pi_k k^z,\tag 2.10$$
and by combining this with (2.7) we will prove
\proclaim{Theorem 2.1} For $z\not=-d+k$, $k=0,1,2,\ldots$, there is
an asymptotic expansion
$$\tr (P_nAP_n)^rQ^z\sim a_r(z)+\sum_{k=d}^{-\infty}a_{r,k}(z)n^{k+z}.
\tag 2.11$$
Moreover, the coefficients in this expansion depend meromorphically on
$z$ and, except for $a_r(z)$, are symbolic functions of $A$. In addition,
$a_{r,k}(z)$ has a simple pole at $z=-k$ and is holomorphic elsewhere,
and $a_r(z)$ is meromorphic with simple poles at $z=-d+k$,
$k=0,1,2,\ldots$.\endproclaim
\demo{Proof} The $j$-th summand above is equal to
$$\tr\sum_{k=1}^{n-\sigma(j)}\bigl(\pi_kA_j\pi_k\bigr)k^z$$
and by (2.7)
$$\tr\sum_{k=1}^{n-\sigma(j)}\bigl(\pi_kA_j\pi_k\bigr)k^z
\sim b(z)+\sum_{l=d}^{-\infty}c_{l-1}\bigl(A_j\bigr)\sum_{k=1}^{n-\sigma(j)}
k^{l-1+z}.$$
By a theorem of Hardy (see [Ha], \S 13.10, page 338)
$$\aligned
\sum_{k=1}^m k^{l-1+z}&\sim C(-z-l+1)+{m^{l+z}-1\over l+z}+{m^{l+z-1}\over 2}\\
&+\sum_{p=1}^\infty (-1)^p(z+l-1)^{(2p-2)}{B_p\over (2p)!}m^{l+z-2p}
\endaligned$$
where $C(s)=\zeta(s)-1/(s-1)$, $B_p$ is the $p$-th Bernoulli
number and $s^{(r)}=s (s+1) \cdots (s+r)$.
Plugging this (with $m=n-\sigma (j)$) into (2.10) we get an expression of the 
form (2.11) where the coefficients are holomorphic in $z$ and $a_{r,k}(z)$ 
is holomorphic except at $z=-k$ where it has a simple pole. Moreover,
if $\hbox{\rm Re} z<-d$ one can take the limit of both sides of (2.11)
as $n$ tends to infinity to obtain
$$\tr A^rQ^z=a_r(z),\tag 2.12$$
and since $\tr A^rQ^z$ is meromorphic with simple poles at $z=-d+k$,
$k=0,1,2,\ldots$, the same is true of $a_r(z)$.
\qed
\enddemo
If we rewrite the right hand side of (2.11) in the form
$$a_r(z)-a_{r,0}(z)+\sum_{k=d,k\not=0}^{-\infty}a_{r,k}(z)n^{k+z}
+za_{r,0}(z){n^z-1\over z}$$
and let $z$ tend to zero we recapture (2.9) with $a_{r,k}=a_{r,k}(0)$
for $k\not=0$, $a_{r,0}=\hbox{\rm Res}_{z=0} a_{r,0}(z)$, and, by (2.12),
$$a_r=\fp\tr A^rQ^z-\fp a_{r,0}(z).\tag 2.13$$
However, $\fp a_{r,0}(z)$ is a local function of $A$ depending only on
the first $d$ terms in its asymptotic expansion; hence the same is
true of $a_r-\fp\tr A^rQ^z$.

Finally by applying this argument to each summand in the series
$$\tr\log (P_nBP_n)Q^z=\sum_{r=1}^\infty {1\over r}\tr 
\bigl(P_nAP_n\bigr)^rQ^z$$
we conclude that the constant term, $b$, in the expansion (1.4) differs
from the zeta regularized ``$\log\det$'' of $B$
$$\fp\tr (\log B)Q^z$$
by a term which is local in $B$ and only depends on the first $d$ terms in 
its symbolic expansion.
\head 3. Zeta regularized determinants\endhead
In this section we relax assumptions on a zeroth order pseudodifferential
operator $B$ and on a regularizer $Q$. We will assume that
the spectrum of $B$ lies in a domain
$D$ of the complex plane where the logarithm is defined and let $\Gamma$ be 
the boundary of $D$ oriented counterclockwise.
Then $\log B$ is defined by the formula
$$\log B={1\over 2\pi i}\int_\Gamma (\lambda I-B)^{-1}d\lambda,$$
and is a zeroth order PDO.
A regularizer $Q$ will be a positive elliptic PDO of order $1$.
The zeta regularized ``$\log\det$'' of $B$ is defined by the formula (1.6).
To compare regularizations of ``$\log\det$'' of $B$ for two different
regularizers, $Q$ and $Q'$, we compute their difference:
$$\align w_Q(B)-w_{Q'}(B)&=\fp \tr\log B\cdot (Q^z-(Q')^z)\\ &=
\res_{z=0}\tr\log B{Q^z-(Q')^z\over z}\\
&=\res[\log B(\log Q-\log Q')].\endalign$$
Notice that $\log B(\log Q-\log Q')$ is a zeroth order pseudodifferential
operator. The last formula shows that $w_Q(B)-w_{Q'}(B)$ is a
local quantity and  depends on the first $d+1$ terms in the symbolic expansions of $B$, $Q$, and
$Q'$.

In the remaining part of this section we will be computing the multiplicative
anomalies for the ``$\log\det$'', namely, $w_Q(AB)-w_Q(BA)$ and
$w_Q(AB)-w_Q(A)-w_Q(B)$. We will show that both are local quantities
and in the case when $d=2$ we will obtain explicit formulas for them
that involve principal symbols of the operators $A$, $B$, and $Q$.
The main tool for computing multiplicative anomalies is the variational
formula for ``$\log\det$''.  
Let $\delta A$ be a variation of an operator and let
$$\sigma_Q(A,\delta A)=\delta w_Q(A)-{1\over 2}\fp\tr\{\delta A A^{-1}+
A^{-1}\delta A\}Q^z.$$
\proclaim{Proposition 3.1}
$\sigma_Q(A,\delta A)$ is a local quantity that depend on $d-1$
terms
in the symbolic expansions of $A$, $\delta A$, and $Q$. If $d=2$ then
$$\sigma_Q(A,\delta A)={1\over 6}
\res (\delta\log a\{\log a,\{\log a,\log q\}\})\tag 3.1$$
where $a(x,\xi)$ is the principal symbol of $A$, $q(x,\xi)$ is the
principal symbol of $Q$, $\{\cdot,\cdot\}$ is the Poisson bracket,
and $\res$ is the symbolic residue (see[Gu].)
\endproclaim
\demo{Proof}
One has
$$\align \sigma_Q(A,\delta A)&={1\over 2\pi i}\int_\Gamma \fp\tr 
  [\log\lambda\Reso \delta A\Reso Q^z]d\lambda\\
  &-{1\over 4\pi i}\int_\Gamma \fp\tr\biggl[{1\over\lambda}\biggl(\delta A\Reso
  +\Reso\delta A\biggr)Q^z\biggr]d\lambda\\
  &={1\over 4\pi i}\int_\Gamma\log\lambda d\lambda \fp\tr 
[2\Reso\delta A\Reso \\ &-
  \delta A(\lambda I-A)^{-2}-(\lambda I-A)^{-2}\delta A]Q^z\\
  &={1\over 4\pi i}\fp\tr\biggl\{\int_\Gamma\log\lambda 
[[\Reso,\delta A],\Reso]d\lambda
  Q^z\biggr\}.\endalign$$
Notice that
$$\tr[[\Reso,\delta A],\Reso]Q^z=\tr [\Reso,\delta A][\Reso, Q^z],$$
and
$$\align\fp &\tr [\Reso,\delta A][\Reso, Q^z]\\
  &=\res_{z=0}\tr {[\Reso,\delta A][\Reso, Q^z]\over z}\\
  &=\res [\Reso,\delta A][\Reso, \log Q].\endalign$$
The operator on the right is of order $d-2$, so its residue
depends on $d-1$ terms in the symbolic expansions of $A$, $\delta A$, and $Q$.
For the variation of ``$\log\det$'' we obtain:
$$\sigma_Q(A,\delta A)=\frac{1}{4\pi i}\int_\Gamma \log\lambda
\res [\Reso,\delta A][\Reso, \log Q]d\lambda.\tag 3.2$$
In the case $d=2$,
$$\align\fp &\tr [\Reso,\delta A][\Reso, Q^z]\\
  &=-\res(\{\reso,\delta a\}\{\reso,\log q\})\\
  &=-\res((\lambda-a)^{-4}\{a,\delta a\}\{a,\log q\}),\endalign$$
and
$$\align\sigma_Q(A,\delta A)&=-{1\over 6}\res (a^{-3}\{a,\delta a\}
\{a,\log q\})
    ={1\over 6}\res(\{a^{-1},\delta a\}\{\log a,\log q\})\\
    &=-{1\over 6}\res(\delta a\{a^{-1},\{\log a,\log q\}\})
    ={1\over 6}\res (\delta\log a\{\log a,\{\log a,\log q\}\}).\endalign$$
\qed
\enddemo
The variation of $w_Q(AB)-w_Q(BA)$ with respect to $A$ (the operator $B$ being 
fixed) equals the sum of 
$$\sigma_Q(AB,\delta AB)-\sigma_Q(BA,B\delta A)\tag 3.3$$
and
$$\align {1\over 2}\fp &(\delta A A^{-1}+B^{-1}A^{-1}(\delta A)B
  -A^{-1}\delta A-B(\delta A)A^{-1}B^{-1})Q^z\\
  &={1\over 2}\res_{z=0}\biggl((\delta A)A^{-1}{Q^z-B^{-1}Q^zB\over z}
-A^{-1}\delta A{Q^z-BQ^z B^{-1}\over z}\biggr)\\
  &={1\over 2}\res((\delta A)A^{-1}(\log Q -B^{-1}\log Q B)
   -A^{-1}\delta A(\log Q-B\log Q B^{-1}))\\
  &={1\over 2}\res ((\delta A)A^{-1}B^{-1}[B,\log Q]+A^{-1}\delta A
[B,\log Q]B^{-1}).\tag 3.4\endalign$$
Both the expressions (3.3) and (3.4) are local and  depend on a finite number of
terms in the symbolic expansions of $A$, $B$, $\delta A$, and $Q$.
Let $A(t)=A^t$. Then the $t$-derivative of $w_Q(A(t)B)-w_Q(BA(t))$
is a local quantity. Clearly, $w_Q(A(0)B)-w_Q(BA(0))=0$; hence
$w_Q(AB)-w_Q(BA)$ is a local quantity. 

The above derivation is valid
if there exists a domain in the complex plane where $\log$ is defined
and that contains the spectrum of $A(t)B$ (and, therefore, of $BA(t)$)
for all $t$, $0\leq t\leq 1$. One can replace the family $A^t$
by any family of zeroth order pseudodifferential operators that connects $A$
with the identity. We will operate under this assumption. It is satisfied if,
for example, the operators $A$ and $B$ are close to the identity or if both of them
are positive.

In the case $d=2$, the quantity (3.3) vanishes because, by (3.1) it depends
on the principal symbols of the operators $A$ and $B$ only, and, on the
level of principal symbols, they commute. Therefore,
$$\align{d\over dt}\bigl( w_Q(A^tB)-w_Q(BA^t)\bigr)&
   ={1\over 2}\res (\log A (B^{-1}[B,\log Q]+[B,\log Q]B^{-1}))\\
   &={1\over 2}\res(\log A(B\log Q B^{-1}-B^{-1}\log Q B)),\endalign$$
and
$$w_Q(AB)-w_Q(BA)={1\over 2}\res(\log A(B\log Q B^{-1}-B^{-1}\log Q B)).
\tag 3.5$$
Now, we fix $A$ and consider the family $B(t)=B^t$. Let
$$g(t)=w_Q(AB^t)-w_Q(B^tA).$$
From (3.5),
$$g'(t)={1\over 2}\res\bigl(\log A(B^t[\log B,\log Q]B^{-t}
       +B^{-t}[\log B,\log Q]B^t)\bigl),$$
and
$$\align g''(t)&={1\over 2}\res \bigl(\log A(B^t[\log B,[\log B,\log Q]]
B^{-t}
       -B^{-t}[\log B,[\log B,log Q]]B^t)\bigl)\\
  &=-{1\over 2}\res (\log a(\{\log b,\{\log b,\log q\}\}-\{\log b,
\{\log b,\log q\}\}))=0;\endalign$$
here $b(x,\xi)$ is the principal symbol of $B$.
In the last equality, we used the fact that the non-abelian residue 
of an operator of order $-2$ on a two-dimensional manifold is the 
residue of the principal symbol.
Hence, 
$$g'(t)=g'(0)=\res (\log A[\log B,\log Q]).$$
Clearly, $g(0)=0$, so
$$w_Q(AB)-w_Q(BA)=\res (\log A[\log B,\log Q]).\tag 3.6$$
The expression (3.6) is of the same form as the Kravchenko--Khesin
cocycle in dimension 1 [KrKh].

The variation of
$$\kappa_Q(A,B)=w_Q(AB)-w_Q(A)-w_Q(B)\tag 3.7$$
with respect to $A$ is the sum of
$$\sigma_Q(AB,\delta A B)-\sigma_Q(A,\delta A)\tag 3.8$$
and
$$\align {1\over 2}\fp&\tr((\delta AB)(AB)^{-1}+(AB)^{-1}\delta (AB)\\
  &-(\delta A) A^{-1}- A^{-1}\delta A)Q^z+\sigma_Q(AB,\delta A B)-
\sigma_Q(A,\delta A)\\
  &={1\over 2}\fp\tr (B^{-1}A^{-1}(\delta A)B-A^{-1}\delta A)Q^z\\
 &={1\over 2}\fp\tr (A^{-1}(\delta A)BQ^zB^{-1}-A^{-1}\delta AQ^z)\\
  &={1\over 2}\res_{z=0}\tr {A^{-1}(\delta A)BQ^zB^{-1}-A^{-1}
\delta AQ^z\over z}\\
   &={1\over 2}\res (A^{-1}(\delta A)B\log Q B^{-1}-A^{-1}\delta A\log Q)\\
  &={1\over 2}\res (A^{-1}\delta A[B,\log Q]B^{-1}).\tag 3.9\endalign$$
Both (3.8) and (3.9) are local expressions and depend on a finite number
of terms in the symbolic expansions of $A$, $B$, $\delta A$, and $Q$.
By taking a family, $A(t)$, that connects $A$ with the identity, we conclude
that $\kappa_Q(A,B)$ is a local quantity ($\kappa_Q(I,B)=0$.)

We will next make these computations more explicit in the two-dimensional situation.
It is covenient to deal with the symmetrized  multiplicative anomaly
$(\kappa_Q(A,B)+\kappa_Q(B,A))/2$. In a similar way to (3.8), (3.9), one
derives
$$\delta\kappa_Q(B,A)=-{1\over 2}\res ((\delta A)A^{-1}B^{-1}[B,\log Q])
+\sigma_Q(BA,B\delta A)-\sigma_Q(A,\delta A),$$
and, therefore,
$$\align \delta{\kappa_Q(A,B)+\kappa_Q(B,A)\over 2}&
  ={1\over 4}\res (A^{-1}\delta A[B,\log Q]B^{-1}-(\delta A)A^{-1}B^{-1}
[B,\log Q])\\
  &+{1\over 2}\sigma_Q(AB, (\delta A)B)+{1\over 2}\sigma_Q(BA, B\delta A)-
\sigma_Q(A,\delta A)\\
 &={1\over 4}\res (\delta A[B,\log Q][B^{-1}, A^{-1}]+\delta A[[B,\log Q], 
A^{-1}B^{-1}])\\
  &+{1\over 2}\sigma_Q(AB, (\delta A)B)+{1\over 2}\sigma_Q(BA, B\delta A)-
\sigma_Q(A,\delta A).\tag 3.10\endalign$$
The first term, $T_1$, on the right in (3.10) equals
$$-{1\over 4}\res (\delta a\{b,\log q\}a^{-2}b^{-2}\{b,a\})
=-{1\over 4}\res (\delta\log a\{\log b,\log q\}\{\log b,\log a\}).$$
The second term, $T_2$, equals
$$\align &-{1\over 4}\res (\delta a\{\{b,\log q\},a^{-1}b^{-1}\})
  =-{1\over 4}\res (\delta\log a\{\{b,\log q\},b^{-1}\})\\
   &\ \ \ \ {1\over 4}\res (b^{-1}\delta a\{\{b,\log q\},a^{-1}\})
  ={1\over 4}\res (b^{-1}\delta\log a\{\{b,\log q\},\log b\})\\
    &+{1\over 4}\res (b^{-1}\delta\log a\{\{b,\log q\},\log a\}).\endalign$$
One uses the identities 
$$\{\{b,\log q\},\log b\}=\{b,\{\log q,\log b\}\}$$
and
$$b^{-1}\{\{b,\log q\},\log a\}=-\{\log a,\{\log b,\log q\}\}+
\{\log b,\log q\}\{\log b,\log a\}$$
to get
$$T_2=-{1\over 4}\res (\delta\log a\{\log (ab),\{\log b,\log q\}\})
+{1\over 4}\res (\delta\log a\{\log b,\log q\}\{\log b,\log a\})$$
and
$$T_1+T_2=-{1\over 4}\res (\delta\log a\{\log (ab),\{\log b,\log q\}\}).$$
By (3.1),
$$\align T_3&={1\over 2}\sigma_Q(AB, (\delta A)B)
 +{1\over 2}\sigma_Q(BA, B\delta A)-\sigma_Q(A,\delta A)\\
&={1\over 6}\res (\delta \log a\{\log (ab),\{\log (ab),\log q\}\})
  -{1\over 6}\res (\delta\log a\{\log a,\{\log a,\log q\}\})\\
&={1\over 6}\res (\delta\log a\{\log a,\{\log b,\log q\}\})
 +{1\over 6}\res (\delta\log a\{\log b,\{\log a,\log q\}\})\\
& +{1\over 6}\res (\delta\log a\{\log b,\{\log b,\log q\}\}).\endalign$$
Finally,
$$\align \delta{\kappa_Q(A,B)+\kappa_Q(B,A)\over 2}&
  =-{1\over 12}\res (\delta\log a\{\log b,\{\log b,\log q\}\})\\
  &-{1\over 12}\res (\delta\log a\{\log a,\{\log b,\log q\}\})\\
  &+{1\over 6}\res (\delta\log a\{\log b,\{\log a,\log q\}\}).\endalign$$
Consider now the family $A(t)=A^t$, the operator, $B$, being fixed. 
Then $\log a(t)=t\log a$,
and
$$\align {d\over dt}&\biggl({\kappa_Q(A(t),B)+\kappa_Q(B,A(t))\over 2}
\biggr)
 =-{1\over 12}\res (\log a\{\log b,\{\log b,\log q\}\})\\
 &-{t\over 12}\res (\log a\{\log a,\{\log b,\log q\}\})
  +{t\over 6}\res (\log a\{\log b,\{\log a,\log q\}\}).\tag 3.11
\endalign$$
The second term on the right in (3.11) vanishes because
$$\res (\log a\{\log a,\{\log b,\log q\}\})
  ={1\over 2}\res (\{\log^2 a,\{\log b,\log q\}\})=0.$$
One integrates (3.11) from $0$ to $1$:
$$\align {\kappa_Q(A,B)+\kappa_Q(B,A)\over 2}&
 ={1\over 12}\res (\log a\{\log b,\{\log(a/b),\log q\}\})\\
&={1\over 12}\res (\{\log a,\log b\}\{\log(a/b),\log q\}\}).\tag 3.12
\endalign$$
(Note that the expression on the right in (3.12) is symmetric in $(a,b)$, as it should be.)

\enddocument